\documentclass[12pt]{amsart}
\usepackage{amsmath}
\usepackage{verbatim}
\usepackage{pagecolor,lipsum}
\usepackage{hyperref}
\usepackage{amsfonts}
\usepackage{graphics}
\usepackage{bbm}
\usepackage{young}
\usepackage{amscd}
\usepackage{color}
\usepackage{amssymb}
\usepackage[all, knot]{xy}
\xyoption{rotate}
\usepackage[x-cp1252]{inputenx}
\usepackage{CJK}
\usepackage{cite}
\usepackage{mathtools}

\xyoption{arc} 

\usepackage{datetime}
\usepackage{amstext}
\usepackage{amsxtra}
\usepackage{amsmath}
\usepackage{amsopn}
\usepackage{amsthm}

\theoremstyle{definition}

\theoremstyle{plain}

\theoremstyle{remark}

\theoremstyle{plain}

\theoremstyle{plain}

\theoremstyle{plain}

\CompileMatrices
\usepackage{listings}
\title{Elementary proofs of generalized continued fraction formulae for $e$}
\author{Zhentao Lu}
\date{\today}
\begin{document}
\begin{abstract}
In this short note we prove two elegant generalized continued fraction formulae 
$$e= 2+\cfrac{1}{1+\cfrac{1}{2+\cfrac{2}{3+\cfrac{3}{4+\ddots}}}}$$
and 
$$e= 3+\cfrac{-1}{4+\cfrac{-2}{5+\cfrac{-3}{6+\cfrac{-4}{7+\ddots}}}}$$ using elementary methods. The first formula is well-known, and the second one is newly-discovered in arXiv:1907.00205 [cs.LG]. We then explore the possibility of automatic verification of such formulae using computer algebra systems (CAS's).
\end{abstract}
\maketitle

\section{Introduction}
We write a (generalized) continued fraction in the form
\begin{equation}
z=b_{0}+{\cfrac {a_{1}}{b_{1}+{\cfrac {a_{2}}{b_{2}+{\cfrac {a_{3}}{b_{3}+{\cfrac {a_{4}}{b_{4}+\ddots \,}}}}}}}},
\end{equation}
with the convergents
\begin{equation}
{\displaystyle z_{0}={\frac {A_{0}}{B_{0}}}=b_{0}, z_{1}={\frac {A_{1}}{B_{1}}}={\frac {b_{1}b_{0}+a_{1}}{b_{1}}}, z_{2}={\frac {A_{2}}{B_{2}}}={\frac {b_{2}(b_{1}b_{0}+a_{1})+a_{2}b_{0}}{b_{2}b_{1}+a_{2}}}, \cdots \,}.
\end{equation}
In general, as stated in \cite{jones1984continued}, $A_n$ and $B_n$ satisfies
\begin{equation}\label{formula:recursive}
\begin{array}{l}
A_{n}=b_{n}A_{n-1}+a_{n}A_{n-2}, B_{n}=b_{n}B_{n-1}+a_{n}B_{n-2}, n\geq 1,\\
{\displaystyle A_{-1}=1, A_{0}=b_{0}, B_{-1}=0, B_{0}=1.} \end{array}
\end{equation}
Thus, there is a straightforward way to verify a proposed continued fraction formula.
Given $\{a_n\},\{b_n\}$, one first computes $A_n$ and $B_n$, and then computes $\displaystyle \lim_{n\to \infty} z_n = \lim_{n\to \infty} \frac{A_n}{B_n}.$

\section{Proofs of two formulae of $e$}
\subsection{The first formula}
We first compute 
\begin{equation}
w= 2+\cfrac{1}{1+\cfrac{1}{2+\cfrac{2}{3+\cfrac{3}{4+\ddots}}}}
= 2+\cfrac{1}{1+\cfrac{\frac{1}{2}}{1+\cfrac{\frac{1}{3}}{1+\cfrac{\frac{1}{4}}{1+\ddots}}}}.
\end{equation}
In this case we have 
\begin{equation}
a_n =  \frac{1}{n}, b_n = 1, n \geq 1.
\end{equation}
Hence
\begin{equation}\label{w:A_nB_n}
\begin{array}{l}
A_{-1} = 1, A_0 = 2, B_{-1} = 0, B_0 = 1,\\
A_n = A_{n-1} + \frac{1}{n}A_{n-2},
B_n = B_{n-1} + \frac{1}{n}B_{n-2}, n\geq 1.
\end{array}
\end{equation}
It is trivial to verify by induction that
\begin{equation}
A_n = n+2 , n\geq 0.
\end{equation}
For $B_n$, we make the auxiliary sequence
\begin{equation}\label{w:k_n}
k_n = \frac{B_{n-2}}{n} - \frac{B_{n-3}}{n-1}, n\geq 2. 
\end{equation}
Then we have 
$k_{n+2} = \frac{1}{(n+2)(n+1)}k_n$ by \eqref{w:A_nB_n}. And in turn this shows $k_n = (-1)^n\frac{1}{n!}$.
Now by \eqref{w:k_n} and \eqref{w:A_nB_n} we have 
\begin{equation}
\sum_{i=2}^{n} k_i = \frac{B_{n-2}}{n},
\end{equation}
hence
\begin{equation}
B_{n} = (n+2) \sum_{i=2}^{n+2} k_i =  (n+2)\sum_{i=2}^{n+2} (-1)^i\frac{1}{i!},
\end{equation}
hence 
\begin{equation}
w= \lim_{n\to \infty}\frac{A_n}{B_n} = \lim _{n\to \infty}\frac{1}{\sum_{i=2}^{n+2} (-1)^i\frac{1}{i!}} = \frac{1}{e^{-1}} = e.
\end{equation}
\subsection{The second formula}
Now we compute
\begin{equation}
v = 3+\cfrac{-1}{4+\cfrac{-2}{5+\cfrac{-3}{6+\cfrac{-4}{7+\ddots}}}}.
\end{equation}
This formula ($v = e$) is recently numerically discovered in \cite{raayoni2019ramanujan} using machine learning techniques.
In this case we have 
\begin{equation}
a_n =  {-n}, \quad b_n = n+3, \quad n \geq 1.
\end{equation}
Hence
\begin{equation}\label{w:A_nB_n}
\begin{array}{l}
A_{-1} = 1, \quad A_0 = 3, \quad B_{-1} = 0, \quad B_0 = 1,\\
A_n = (n+3)A_{n-1} - {n}A_{n-2},\quad 
B_n = (n+3)B_{n-1} - {n}B_{n-2},\quad  n\geq 1.
\end{array}
\end{equation}
This time it is trivial to verify by induction that
$B_{n} = \frac{(n+1)^2}{n}B_{n-1}$, so
\begin{equation}
B_n = \frac{(n+1)!^2}{n!} = (n+1)\cdot (n+1)!, \quad n\geq 1.
\end{equation}
For $A_n$, an observation of the first few terms suggests that it is a shifted version of the sequence A001339 in \cite{sloane2003line}. Then it is routine to verify that 
\begin{equation}\label{v:A_n-induction}
A_n = \sum_{k=0}^{n+1} (k+1)!\begin{pmatrix}
n+1\\k
\end{pmatrix}, \quad n\geq 1,
\end{equation}
by induction argument.
\footnote{We find that the inductive step is easier to carry out if using the equivalent expression 
\begin{equation}
A_n = \sum_{k=0}^{n+1} (n+2-k)!\begin{pmatrix}
n+1\\n+1-k
\end{pmatrix}.
\end{equation}
}

With the expression of $A_n$ and $B_n$ known, we get
\begin{equation}
v_n =  \frac{A_n}{B_n} 
=\sum_{k=0}^{n+1} \frac{(k+1)}{(n+1)\cdot(n+1-k)!}
=\sum_{k=0}^{n+1} \frac{(n+2-k)}{(n+1)\cdot k!}.
\end{equation}
Hence
\begin{equation}
v=\lim_{n\to\infty} v_n 
= \lim_{n\to\infty} \frac{n+2}{n+1}\sum_{k=0}^{n+1}\frac{1}{k!}
  - \frac{1}{n+1}\sum_{k=0}^{n}\frac{1}{k!} 
= e.
\end{equation}
\section{Computer-aided verification of generalized continued fraction formulae}
Given a generalized continued fraction formula, we propose the following work flow to verify it:

\begin{tabular}{l}
Step 1. extract the $\{a_n\},\{b_n\}$ terms.\\
Step 2. Get the recursive formulae of $A_n$ and $B_n$ using \eqref{formula:recursive}.\\
Step 3. Compute the first few terms of $A_n$ and $B_n$ and using \\
\quad \quad \quad \quad WolframAlpha and OEIS to guess a closed-form of them.\\
Step 4. Using mathematical induction to prove the closed-form\\ 
\quad \quad \quad \quad expression.\\
Step 5. Compute the limit and check against the proposed formula.
\end{tabular}

\section*{}
{\bf Contact:} Zhentao Lu, zhentao@sas.upenn.edu

\bibliographystyle{amsplain}

\end{document}